\documentclass[12pt,a4paper]{article}
\usepackage[cp850]{inputenc}
\usepackage{latexsym,amsfonts}

\newtheorem{Thm}   {Theorem}  [section]
\newtheorem{Lem}[Thm]{Lemma}

\newtheorem{Rem}[Thm]{Remark}

\newtheorem{Def}[Thm]{Definition}
\newtheorem{tab}[Thm]{Table}

\author{ Vladimir D. Tonchev\thanks{ Research supported by NSA Grant
H98230-16-1-0011.}\\
 Michigan Technological University\\
 Houghton, Michigan 49931, USA
}

\title{\bf Linearly Embeddable Designs}

\begin{document}
\maketitle
\begin{abstract}

A residual design ${\cal{D}}_B$ with respect to a block $B$ of a
given design $\cal{D}$ is defined to be linearly embeddable
over $GF(p)$ if the $p$-ranks of the incidence matrices
of  ${\cal{D}}_B$ and $\cal{D}$ differ by one.
A sufficient condition for a residual design to be linearly embeddable is proved
in terms of the minimum distance of the linear code spanned by the
incidence matrix, and this condition is used to show
that the residual designs of several known infinite classes of designs
are linearly embeddable.
A necessary condition for linear embeddability is proved for affine resolvable
designs and their residual designs.
As an application, it is shown that a residual design of the classical affine design
of the planes in $AG(3,2^2)$ admits two nonisomorphic embeddings over $GF(2)$
that give rise to the only known counter-examples of Hamada's conjecture
over a field of non-prime order.

{\bf Keywords:} incidence matrix, residual design, affine resolvable design, 
linear code, $p$-rank.

{\bf Mathematics Subject Classification:} 05B05, 05B25, 51E20, 94B05. 

\end{abstract}

\vspace{5mm}

\section{Introduction}

We assume familiarity with basic facts and notions from
combinatorial design theory and coding theory
 (\cite{AK}, \cite{BJL}, \cite{BM}, \cite{Hall}, \cite{HP}, \cite{Ton-HCT},
\cite{Ton-CRC},  \cite{T88}).

A combinatorial {\it design} (or an {\it incidence structure}) is a pair
$\cal{D}$=$(X, \cal{B})$ of a finite set of {\it points} $X=\{ x_i \}_{i=1}^v$
and a collection $\cal{B}$=$\{ B_j \}_{j=1}^b$ of subsets $B_j \subseteq X$, called
{\it blocks}. The point-by-block incidence matrix $A=( a_{i,j})$ of $\cal{D}$=$(X, \cal{B})$
is a $v$ by $b$ $(0,1)$-matrix with
$a_{i,j}=1$ if $x_i \in B_j$, and $a_{i,j}=0$ otherwise.

An incidence structure is {\it simple} if $X$ is a proper set consisting of $v$ distinct points,
and all blocks are distinct subsets of points, or equivalently, its incidence matrix does not have 
any identical rows or columns.

Given integers $v\ge k \ge t \ge 0$, $\lambda\ge 0$, a $t$-$(v,k,\lambda)$ design (or briefly, a $t$-design)
 $\cal{D}$
is an incidence 
structure with $v$ points and blocks of size $k$ such that every $t$-subset of points is contained
in exactly $\lambda$ blocks. It follows that a $t$-$(v,k,\lambda)$ is also an $s$-$(v,k,\lambda_s)$
design for any $0 \le s \le t$, where
\[ \lambda_s =\frac{ { v-s \choose t-s }}{{k-s \choose t-s}}. \]
In particular, the total number of blocks is given by
\[ b=\lambda_0 = \frac{ { v \choose t }}{{k \choose t}}. \]
The number $\lambda_1$ of blocks containing a given point is often called the {\it replication number},
and is denoted by $r$.

The number of blocks $b$ and the number of points $v$ of any 2-$(v,k,\lambda)$ 
design $D$ with $v>k>0$ satisfy the following
inequality, known as the Fisher inequality:
\begin{equation}
\label{Fish}
 b \ge v, \end{equation}
where the equality $b=v$ holds if and only if every two distinct blocks of $D$
share exactly $\lambda$ points.
A 2-$(v,k,\lambda)$ design $D$ with $b=v$ is called {\it symmetric}.
If $A$ is the incidence matrix of a symmetric 2-$(v,k,\lambda)$ design $D$, 
then $A^{T}$ is the incidence matrix of a symmetric design with the same parameters, called
the {\it dual} of $D$.

%\section{Derived and residual designs}

Let $\cal{D}$=$(X,\cal{B})$ be design, and let $B \in \cal{B}$ be a block of $\cal{D}$.
The incidence structure
\[ {\cal{D}}^B = (X',\cal{B'}), \]
where
\[ X'=B, \ {\cal{B'}}= \{ B\cap B_j \ | \ B_j \in {\cal{B}}, B_j \neq B \}, \]
is called the {\it derived} design of $\cal{D}$ with respect to block $B$.

Accordingly, the incidence structure 
\[ {\cal{D}}_B=(X'', \cal{B''}), \]
 where
\[ X'' = X\setminus B, \ {\cal{B''}}=\{ B_j \setminus (B_j\cap B) \ | \  B_j \in {\cal{B}}, B_j \neq B \}, \]
is called the  {\it residual} design of $\cal{D}$ with respect to block $B$.

The main subject of this paper are incidence structures $\cal{D}$ with the property
that for some prime number $p$, the $p$-ranks of the incidence matrices of $\cal{D}$
and some of its residual designs, ${\cal{D}}_B$, differ by one.
We prove a sufficient condition, as well as some necessary conditions
for an incidence structure to admit this property. As an application, we give
an alternative construction of the 2-$(64,16,5)$ designs of 2-rank 16,
being the only known counter-examples of Hamada's conjecture
\cite{Ha}, \cite{Ha81} over a field 
of non-prime order \cite{CJT}, \cite{HLT}, \cite{JT09}, \cite{T86}.

\section{Linearly embeddable residual designs}
\label{lin}

Suppose that $\cal{D}$=$(X,\cal{B})$ is a design, 
and $B\in{\cal{B}}$ is a block containing
$k$ points, $k > 1$. For convenience of notation, we assume that the points and blocks
of $\cal{D}$ are labeled so that $B$ is the last block, and consists of 
the first $k$ points of $X$.
Then the point-by-block $v \times b$ incidence matrix $A$ of $\cal{D}$ 
can be written as in eq. (\ref{eq1}), where $A'$ is the $k \times (b-1)$ 
incidence matrix of the derived
design ${\cal{D}}^B$, and
$A''$ is the incidence matrix of the residual design ${\cal{D}}_B$.

\begin{equation}
\label{eq1}
A=\left( \begin{array}{ll}
 & 1 \\
& \cdot \\
 A' & \cdot \\
& \cdot \\
 & 1 \\
 & 0\\
& \cdot \\
A'' & \cdot \\
& \cdot  \\
 & 0\\
\end{array}\right),
\end{equation}

It is obvious from (\ref{eq1}) that
\begin{equation}
\label{eq2}
rank_p A \ge rank_p A'' +1,
\end{equation}
where $rank_p$ denotes the $p$-rank of the corresponding matrix,
that is, its rank
over a finite field $GF(p)$ of prime order $p$.

\begin{Def}
{\rm
We say that a residual design ${\cal{D}}_B$ is {\it linearly embeddable} 
over $GF(p)$ if
\begin{equation}
\label{eq3}
rank_p A = rank_p A'' +1.
\end{equation}
}
\end{Def}

The condition (\ref{eq3}) implies that all rows of $A$ belong to the
linear code of length $b$ over $GF(p)$, being the row span
of a $(v-k+1)\times b$ matrix, consisting of the $v-k$ rows of $A$
labeled by the points of ${\cal{D}}_B$, plus one extra row 
$y=(y_1,\ldots, y_b)$ from 
the row span of $A$ such that $y_b\neq 0$ (for example, $y$ can be
any of the $k$ rows of $A$ labeled by a point of ${\cal{D}}^B$).

Clearly, the condition (\ref{eq3}) is a strong requirement. For example,
this condition does not hold if $rank_p A = v$ and $k\ge 2$.

The next theorem gives a sufficient condition for a residual design 
to be linearly embeddable.

\begin{Thm}
\label{t1}
Let $\cal{D}$=$(X,{\cal{B}})$ be a design with $v$ points, $b$ blocks,
and a $v \times b$ incidence matrix $A$, and let $C$ be the linear code
of length $v$ over $GF(p)$ spanned by the columns of $A$.
If the minimum Hamming weight of $C$ is $d$, then every residual design
${\cal{D}}_B$ with respect to a block of size $d$ is linearly
embeddable over $GF(p)$.
\end{Thm}

{\bf Proof}. 
Let $y\in C$ be a codeword of minimum weight $d$, such that the support of $y$
(that is, the set of indices of its nonzero components)
is a block $B$ of $\cal{D}$. We consider the linear code $C_y$ of length $v-d$
obtained from $C$ by puncturing the $d$  coordinates labeled by the support
of $y$, or equivalently, by the points of $B$.
Clearly, the dimension of $C_y$ is equal to $rank_p A''$, where $A''$
is the incidence matrix of the residual design ${\cal{D}}_B$,
while the dimension of $C$ is equal to $rank_p A$.

In the terminology \cite{HS}, 
and the notation of \cite[Section 2.7., page 80]{HP}, $C_y$ is
the residual code $Res(C,y)$ with respect to $y$.
Since 
\[ d < \frac{p}{p-1}d, \]
 it follows from a result by Hill and Newton \cite[Lemma 2.13]{HN}
(see also \cite[Theorem 2.7.1, page 80]{HP}) that the dimension of $C_y$ is 
smaller than the dimension of $C$ by one, hence
\[ rank_p A'' = rank_p A -1, \]
which completes the proof.
$\Box$

Next we give examples of 2-designs which satisfy the condition of
Theorem \ref{t1}, and consequently, have linearly embeddable residual designs.

A symmetric 2-design with parameters
\begin{equation}
\label{sSDP}
 2-(2^{2m}, 2^{2m-1} - 2^{m-1}, 2^{2m-2} - 2^{m-1}) 
\end{equation}
has the {\it symmetric difference property}, and is called an SDP design 
(Kantor \cite{K})
if the symmetric difference of any three blocks is either a block or a complement
of a block. The number of nonisomorphic SDP designs with
parameters (\ref{sSDP}) grows exponentially with $m$
\cite{K2}.

Dillon and Schatz \cite{DS} proved the following characterization of symmetric SDP designs:
a design with parameters (\ref{sSDP}) has the symmetric difference property if and only if
its blocks are the supports of minimum weight codewords in a binary linear code
of length $2^{2m}$ and dimension $2m+2$, spanned by a bent 
function on $2m$ variables and the first order Reed-Muller code of length $2^{2m}$.
This result and Theorem \ref{t1} imply the following.

\begin{Thm}
\label{t2}
The residual designs of a symmetric SDP design are linearly embeddable over $GF(2)$.
\end{Thm}

A residual design of a symmetric SDP design $D$ with parameters (\ref{sSDP}) is a 2-design
with parameters
\begin{equation}
\label{sdp-r}
2-(2^{2m-1}+2^{m-1}, 2^{2m-2}, 2^{2m-2} - 2^{m-1}), \end{equation}
and derived design of $D$ has parameters
\begin{equation}
\label{sdp-d}
2-(2^{2m-1}-2^{m-1}, 2^{2m-2} - 2^{m-1}, 2^{2m-2} - 2^{m-1}-1). 
\end{equation}
The residual and derived designs of a symmetric SDP design have the property that
the symmetric difference of any two blocks is either a block or a complement of block,
and are called quasi-symmetric SDP designs (the term "quasi-symmetric" refers to the property
that there are only two distinct block intersection numbers; see Shrikhande \cite{MSS} for an 
introduction to quasi-symmetric designs).
It was shown by Jungnickel and the author in \cite{JT} that the number
of nonisomorphic quasi-symmetric SDP designs with parameters (\ref{sdp-r}) or (\ref{sdp-d})
grows exponentially with $m$.

\begin{Rem}
{\rm
It was proved by the author in \cite{Tongd} that any quasi-symmetric SDP design is
uniquely embeddable as a residual or derived design in a symmetric SDP design.
Combined with the result of Theorem \ref{t2}, this implies that a residual design
of a symmetric SDP design $D$ is linearly embeddable in a unique (up to isomorphism)
symmetric design, namely $D$. We will discuss some interesting linearly embeddable
residual designs later in this paper, which can be linearly embedded in
two nonisomorphic designs.
}
\end{Rem}

By the Dillon-Schatz theorem \cite{DS}, 
the 2-rank of a symmetric SDP design
with parameters (\ref{sSDP}) is $2m+2$, and consequently,
the 2-rank of its residual or derived designs is equal to $2m+1$.

The binary linear code $C'$  of length $2^{2m-1}-2^{m-1}$ spanned by the blocks of a quasi-symmetric
SDP design with parameters (\ref{sdp-d}) consists of the zero vector, the all-one vector,
the incidence vectors of the blocks (which are also the minimum weight codewords), and the
incidence vectors of the complements of the blocks.
Similarly, the binary linear code $C''$ of length $2^{2m-1}+2^{m-1}$ spanned by the blocks
of a quasi-symmetric
SDP design with parameters (\ref{sdp-r}) consists of the zero vector, the all-one vector,
the incidence vectors of the blocks (which are also the minimum weight codewords), and the
incidence vectors of the complements of the blocks.
Thus, the Dillon-Schatz theorem \cite{DS} and Theorem \ref{t1} imply the following.

\begin{Thm}
\label{t3}
The residual design with respect to any block of a given quasi-symmetric SDP design
is linearly embeddable over $GF(2)$.
\end{Thm}

Other infinite classes of linearly embeddable designs over $GF(2)$ are
the designs supported by the minimum weight codewords of  Reed-Muller
codes, or punctured Reed-Muller codes.

The codewords of minimum weight $d=2^{m-r}$ of the Reed-Muller code of length $2^m$
and order $r$ ($1 \le r <m-1$), form the block by point incidence matrix of a 3-design $\cal{D}$,
which is linearly embeddable over $GF(2)$ by Theorem \ref{t1}.
In addition, puncturing one of the $2^m$ coordinates of the Reed-Muller code of length $2^m$
and order $r$ gives a punctured code $C'$ of length $2^m -1$ and minimum distance $2^{m-r} -1$.
The minimum weight codewords of $C'$ support a linearly embeddable 2-design by Theorem \ref{t1}.

The linearly embeddable designs supported by the binary Reed-Muller
codes and punctured Reed-Muller
codes are special classes of designs based on finite geometry.
Let $q=p^t$, $p$ - prime, $t\ge 1$.
The $d$-dimensional subspaces of the $n$-dimensional
projective geometry $PG(n,q)$ over $GF(q)$, are the blocks of a 2-design, denoted by $PG_{d}(n,q)$,
 with parameters
\[2-\left(\frac{q^{n+1}-1}{q-1}, \frac{q^{d+1}-1}{q-1}, \frac{(q^{n+1}-q^{2})(q^{n+1}-q^{3})\cdots
(q^{n+1}-q^{d})}{(q^{d+1}-q^{2})(q^{d+1}-q^{3})\cdots(q^{d+1}-q^{d})} \right)\]

Similarly, the $d$-subspaces of the $n$-dimensional affine geometry $AG(n,q)$ over $GF(q)$
are the blocks of a 2-design, denoted by$AG_{d}(n,q)$, with parameters

\[2-\left(q^{n}, q^{d}, \frac{(q^{n}-q)(q^{n}-q^{2})\cdots(q^{n}-q^{d-1})}{(q^{d}-q)(q^{d}-q^{2})
\cdots(q^{d}-q^{d-1})} \right)\]

If $q=2$,  $AG_{d}(n,2)$ is also a $3$-$(2^n, 2^d, \lambda_3)$ design with
\[ \lambda_3=\frac{(2^n -2^2)\cdots (2^n -2^{d-1})}{(2^d -2^2)\cdots (2^d -2^{d-1})}. \]  

The incidence vectors of the blocks of $PG_{d}(n,q)$ are codewords of minimum
weight of the $GF(p)$-subfield subcode of a nonprimitive generalized Reed-Muller
code over $GF(q)$ \cite[5.7]{AK}, \cite[Chapter 2]{BM}.
Similarly, the blocks of $AG_{d}(n,q)$ are supported by minimum weight codewords
in the code over $GF(p)$ spanned by the incidence vectors of the
blocks \cite[Corollary 5.5.5, page 166]{AK}. Thus, by Theorem \ref{t1},
 we have
\begin{Thm}
\label{t4}
The residual designs of $PG_{d}(n,p^t)$ and $AG_d(n,p^t)$ are linearly embeddable over $GF(p)$.
\end{Thm}

\section{Residual designs of affine resolvable designs}

The Fisher inequality (\ref{Fish}) for a 2-$(v,k,\lambda)$ design with $v>k>0$
can be strengthened when $v$ is a multiple of $k$, $v=qk$, as follows:
\begin{equation}
\label{resol}
b \ge v+r -1, \end{equation}
where $r=(v-1)\lambda/(k-1)$ is the replication number.
Suppose that $\cal{D}$ is a 2-$(qk,k,\lambda)$ design with $q>1$ and $k>0$. Any set of $q$ pairwise
disjoint blocks is called a {\it parallel class}. A
{\it resolution} of $\cal{D}$ is a partition of the collection of blocks into
$r$ disjoint parallel classes. A design is {\it resolvable} if it admits at least one resolution.

The parameters of a resolvable 2-design satisfy the inequality (\ref{resol}).
In addition, a 2-$(qk,k,\lambda)$ design with $b=v+r-1$ is resolvable
if and only if the number
\[ \mu = \frac{k}{q} = \frac{k^2}{v} \]
is an integer, and every two blocks are either disjoint or share exactly $\mu$ points
(cf. Bose \cite[Theorem 1.6.1]{Bose}, or \cite[Theorem 2.3.3]{T88}).

A resolvable 2-$(qk,k,\lambda)$ design with $b=v+r-1$ blocks is called {\it affine resolvable}.
An affine resolvable design admits only one resolution, and its parameters can be written as
\begin{equation}
\label{afres}
v=q^2\mu, \ k=q\mu, \ \lambda=\frac{q\mu -1}{q-1}.
\end{equation}
If $B$ is a block of an affine resolvable
2-$(q^2\mu, q\mu, \frac{q\mu -1}{q-1})$ design $\cal{D}$, the derived design ${\cal{D}}^B$
is a 2-$(q\mu, \mu, \frac{q\mu -1}{q-1} - 1)$ design (here we do not consider the empty intersections
of the $q-1$ blocks from the parallel class of $B$ as blocks of the derived design ${\cal{D}}^B$).

 Any affine geometry design $AG_{d}(n,q)$, $1\le d \le n-1$, is resolvable: 
one resolution has as parallel classes
the collections of cosets of the affine $d$-subspaces through the origin.

The number of nonisomorphic resolvable designs having the same
parameters as $A_d(n,q)$, $3\le d\le n-1$, grows exponentially 
 (Jungnickel \cite{J84}, Lam, Lam and Tonchev \cite{LLT}).

 If $d=n-1\ge 1$, 
 $A_{n-1}(d,q)$ is an affine resolvable 2-$(q^n, q^{n-1}, \frac{q^{n-1}-1}{q-1})$ design.
If $B$ is a block of $\cal{D}$=$AG_{n-1}(n,q)$, there are $q-1$ blocks parallel to $B$,
and every other block intersects $B$ in a $(n-2)$-subspace of $AG(n,q)$. The non-empty
intersections of $B$ with other blocks
 of $\cal{D}$ form a 2-design ${\cal{D}}'$ with point set $B$
and parameters
\[ v' = q^{n-1}, \ k' = q^{n-2}, \ \lambda' =\frac{q^{n-1}-1}{q-1}-1 =  
 \frac{q(q^{n-2}-1)}{q-1}, \ b'=\frac{q^{2}(q^{n-1}-1)}{q-1}. \]
The design  ${\cal{D}}'$ is not simple: its collection of blocks
is a multi-set, where every block appears with multiplicity $q$.
A set $S$ of distinct representatives of the blocks of  ${\cal{D}}'$
consists of all $(n-2)$-subspaces of $AG(n-1,q)$, where the points
of $AG(n-1,q)$ are identified with the points of $B$. Thus, $S$ is a simple
2-$(q^{n-1},q^{n-2},\frac{q^{n-2}-1}{q-1})$ design isomorphic to
$AG_{n-2}(n-1,q)$.

The residual design ${\cal{D}}_B$ has $q-1$ blocks of size $q^{n-1}$ (these are  blocks
from the parallel class of $\cal{D}$ that contains $B$), while
the remaining
\[ \frac{q(q^n -1)}{(q-1)}-q \]
blocks of ${\cal{D}}_B$ are of size $q^{n-1} - q^{n-2}$.

Let  ${\cal{D}}''$ be the substructure of ${\cal{D}}_B$
consisting of all blocks of size $q^{n-1} - q^{n-2}$.
Since every two non-parallel blocks of  $\cal{D}$ meet in $\mu=q^{n-2}$ points,
each set of $q$ identical blocks of  ${\cal{D}}'$ corresponds to a parallel class of $q$
pairwise disjoint blocks of  ${\cal{D}}''$. In this way, we obtain a resolution $R$ of
 ${\cal{D}}''$, in which the parallel classes are labeled by the blocks of $S$.

This construction can be applied to residual designs of other affine resolvable
designs having the parameters of $AG_{n-1}(n,d)$, provided that there is a block
satisfying the condition of the following definition.

\begin{Def}
\label{good}
{\rm
A block $B$ of an affine resolvable 2-$(q^n, q^{n-1}, \frac{q^{n-1} -1}{q-1})$ design
$\cal{D}$ is called a {\it good} block if the nonempty intersections of $B$ with the 
remaining blocks form
a 2-$(q^{n-1}, q^{n-2}, \frac{q(q^{n-2} -1)}{q-1})$ design ${\cal{D}}'$, whose
collection of blocks is a union of $q$ identical copies of the block set of a simple 
2-$(q^{n-1}, q^{n-2}, \frac{q^{n-2} -1}{q-1})$ design $S$.  
}
\end{Def}

\begin{Rem}
{\rm
The above definition of a good block concerns a special case of a more general concept
introduced by Kimberley \cite{Kim} and used by Kantor \cite{K69} (see Beth,  Jungnickel 
and Lenz \cite[XII.5]{BJL}) for further references).
}
\end{Rem}

We note that
by the inequality of Mann \cite{Mann},
\cite[Theorem 1.1.5, page 6]{T88},
every 2-$(q^{n-1},q^{n-2},\frac{q^{n-2}-1}{q-1})$ design is simple.

Clearly, any good block $B$ defines 
a resolution  of the subdesign ${\cal{D}}''$ of the residual design
${\cal{D}}_B$, consisting of the blocks of size $q^{n-1} - q^{n-2}$.

\begin{Thm}
\label{t5}

Let $\cal{D}$ be an affine resolvable 2-$(q^n, q^{n-1}, (q^{n-1}-1)/(q-1)$)
design, ($n \ge 2$), with a good block $B$, where $q=p^t$, $p$ is prime, and  $q\ge 4$.
If the residual design ${\cal{D}}_B$ is linearly embeddable over $GF(p)$, then
the linear code over $GF(p)$ of length $q^{2}(q^{n-1}-1)/(q-1)$, spanned by the rows 
of the incidence matrix $M''$
of the substructure ${\cal{D}}''$ of  ${\cal{D}}_B$ consisting of all blocks of size 
$q^{n-1} - q^{n-2}$,
contains at least $(p-1){ q^{n-1} \choose 2}$
codewords of weight $2q^{n-1}$, whose supports are unions of parallel classes
of the resolution $R$ of  ${\cal{D}}''$ defined by $B$.
 
\end{Thm}

{\bf Proof.}
For convenience of notation, we assume that the
points and blocks of $\cal{D}$ are labeled so that $B$ is the last block and
consists of the first $q^{n-1}$ points. Then the point-by-block incidence 
matrix $A$ of $\cal{D}$ is given by (\ref{eq1}), where $A'$ is the
incidence matrix of the derived design ${\cal{D}}^B$, and
$A''$ is the incidence matrix of the residual design ${\cal{D}}_B$. 
Let
\begin{equation}
\label{a1a2}
 A_1=\left( \begin{array}{ll}
 & 1 \\
& \cdot \\
 A' & \cdot \\
& \cdot \\
 & 1 \\
\end{array}\right), \ A_2 = \left( \begin{array}{ll} 
 & 0\\
& \cdot \\
A'' & \cdot \\
& \cdot  \\
 & 0\\
\end{array}\right).
\end{equation}
The matrix $A'$ contains $q-1$ all-zero columns that correspond to 
the blocks of $\cal{D}$ parallel to $B$.
We denote by $M'$ the submatrix of all nonzero columns of $A'$.
 
Since ${\cal{D}}_B$ is linearly embeddable over $GF(p)$,  we have
\[ rank_{p}A_2 =rank_{p}A -1, \]
which implies that the vector space $L_2$ being the span of the rows of $A_2$
over $GF(p)$, coincides with the subspace of co-dimension 1 of the row span $L$ of $A$, 
consisting of all vectors in $L$ having last coordinate equal to zero.
It follows that the difference of any two rows of $A_1$ belongs to $L_2$, 
and consequently, the difference
of any two rows of $M'$ belongs to the rows space of $M''$.

The set of the
\[ r' =q\frac{q^{n-1}-1}{q-1} \]
nonzero positions of any row of $M'$ corresponds to a union of $(q^{n-1}-1)/(q-1)$ parallel classes 
of the resolution $R$. 
Every two distinct rows of $M'$ overlap in a set $T$ of
\[ \lambda' = q\frac{q^{n-2}-1}{q-1} \]
nonzero positions, where $T$ corresponds to a union of $(q^{n-2}-1)/(q-1)$ parallel classes of $R$.
Thus, the difference of every two distinct rows of $M'$ is a vector of Hamming weight
\[ 2(r' -\lambda')=2q^{n-1}, \]
belonging to the row space of $M''$, whose support is a union of $2q^{n-2}$ parallel classes of $R$.
 
We will show that the differences of different pairs of distinct rows of $A'$ are distinct vectors
of weight $2q^{n-2}$. Let ${\cal{P}}_1$ =$\{ r_1, r_2 \}$, ($r_1 \neq r_2$),  
${\cal{P}}_2$ =$\{ r_3, r_4 \}$,  ($r_3 \neq r_4$), be two distinct unordered pairs of rows.
Fist, suppose that  ${\cal{P}}_1$ and ${\cal{P}}_2$ share one row. If $r_1 = r_3$ then
\begin{equation}
\label{dif}
 r_1 - r_2 = r_3 - r_4
\end{equation} 
implies $r_2 = r_4$ and  ${\cal{P}}_1$=${\cal{P}}_2$, a contradiction.
If $r_1 = r_4$ then $r_2 \neq r_3$, and equation (\ref{dif}) implies $2r_1 =r_2 + r_3$, 
   which is impossible due to the Hamming weights of the rows and the size
 of the overlap of their supports.
   Hence, if the pairs ${\cal{P}}_1$,  ${\cal{P}}_2$ comprise of three distinct rows of $A'$, we have
\begin{equation}
\label{dif2}
r_1 - r_2 \neq r_3 - r_4.
\end{equation}
Suppose now that $r_1, r_2, r_3, r_4$ are four distinct rows of $A'$ which satisfy
the equation (\ref{dif}). Let $C$ be the linear code of length $q^{n-1}$ over $GF(p)$,
being the null space of the column space of $A'$. The equation (\ref{dif}) implies
that $C$ contains a codeword of weight 4 with support labeled by the four rows, 
thus, the minimum weight of $C$ is at most 4. 
 We will show, however, that if $q\ge 4$, the minimum weight $d$ of $C$ is at least 5.
 
Let  ${M'}_S$ be a submatrix of $M'$ consisting of $q(q^{n-1}-1)/(q-1)$ distinct 
columns of $M'$, that is,  ${M'}_S$ is an incidence matrix of the simple
2-$(q^{n-1}, q^{n-2}, \frac{q^{n-2}-1}{q-1})$ subdesign $S$ of ${\cal{D}}'$.

By Rudolph's theorem \cite{Rud}, 
\cite[Theorem 2.7.3]{T88},
the code $C$ can correct up to
\[ e=\lfloor \frac{r+\lambda -1}{2\lambda}\rfloor \]
errors by a majority-logic decoding algorithm using the columns of ${M'}_S$,
where
\[ r = \frac{q^{n-1} -1}{q-1} \]
is the replication number of $S$, and
\[ \lambda =  \frac{q^{n-2}-1}{q-1} \]
is the number of blocks through a pair of points.

We have
\[ \frac{r+\lambda -1}{2\lambda} \ge \frac{r}{2\lambda} =\frac{q^{n-1} -1}{2(q^{n-2}-1)} \ge 2, \]
provided that $q\ge 4$ and $n\ge 2$.
Thus, the code $C$ can correct at least 2 errors, which implies that
the weight of any nonzero codeword of $C$ is greater than or equal to 5.

Hence, the differences of pairs of rows of $A'$ are all different
codewords from the row space
of $M''$, each of weight $2q^{n-1}$, and having a support being a union of parallel classes of $R$.
Taking into account the $p-1$ nonzero scalar multiples of each such codeword gives 
a set of $(p-1){ q^{n-1} \choose 2}$ distinct codewords with the required property.
This completes the proof. {$\Box$}

\section{Residual designs of $AG_{2}(3,4)$}
\label{sec5}

The smallest parameters $n$, $q$ that satisfy the conditions of Theorem \ref{t5}
are $n=2$ and $q=4$. Any residual design of the 2-$(64,16,5)$ design
 $\cal{D}$=$AG_{2}(3,4)$ is linearly embeddable
over $GF(2)$ by Theorem \ref{t4}. The 84 blocks of $\cal{D}$ are the planes
in $AG(3,4)$, and all blocks are in one orbit under the collineation group 
of $AG_{2}(3,4)$, being of order
\[ 2\cdot 4^3(4^3-1)(4^3 -4)(4^3 -4^2) = 23,224,320 = 2^{13}\cdot 3^4 \cdot 5 \cdot 7. \]
Thus, all residual designs of $\cal{D}$ are isomorphic.
By Hamada's formula \cite{Ha}, the 2-rank of $AG_{2}(3,4)$
is 16.

Let $B$ be a block of $\cal{D}$.
Our goal is to determine if ${\cal{D}}_B$ =$(AG_{2}(3,4))_B$ can be embedded linearly
over $GF(2)$ as a residual design with respect to a good block into any other affine
resolvable 2-$(64,16,5)$ design $\cal{E}$ which is not isomorphic to
$AG_{2}(3,4)$.

 The weight distribution 
of the binary linear code of length 80 and dimension 15
being the row space of the $48 \times 80$ incidence matrix
$M''$ of the substructure ${\cal{D}}''$ of the residual design 
${\cal{D}}_B$ consisting of all 
blocks of size 12, is given in Table \ref{tab1}.

\begin{tab}
\label{tab1}
\end{tab}
\begin{tabular}{|r|r|r|r|r|r|r|r|r|r|r|r|r|}
\hline
$i$ & 20 & 30 & 32 & 34 & 36 & 38 & 40 & \ldots & 48 & 50 & 52 & 64 \\
\hline
$A_i$ & 48 &  768 & 610 &  1280 & 6240 & 7680 &  2880 & \ldots & 600 & 256 & 240 & 5\\
\hline
\end{tabular}

\vspace{5mm}

The design ${\cal{D}}''$ has 40 parallel classes and 32 resolutions.
The automorphism group of ${\cal{D}}''$ is of order $552,960$, and partitions the set of 32
resolutions in three orbits of lengths 2, 10, and 20, respectively.
The resolution $R$ whose parallel classes are labeled by the blocks of the
2-$(16,4,1)$ subdesign of ${\cal{D}}^B$ is one of the resolutions in the orbit of length 2, 
the second one being the resolution induced by the unique resolution of $\cal{D}$.
 
It is easy to verify that among the 610 codewords of weight 32 (cf. Table \ref{tab1}),
there are 130 codewords
whose supports are unions of 8 parallel classes from  a resolution
from the orbit of length 2, 34 such codewords with respect to a
resolution from the orbit of length 10, and 10 codewords with respect to a resolution
from the orbit of length 20. By Theorem \ref{t5}, this implies that one can have
a linear embedding only with respect to a resolution from the orbit of length 2.
Thus, it is sufficient to consider linear embeddings with respect to the resolution $R$.

To search for such linear embeddings, we extend the $48 \times 80$ incidence matrix of
${\cal{D}}''$ by four columns: three columns of weight 16, being the incidence vectors
of the three blocks of $\cal{D}$ parallel to $B$, plus one all-zero column.
Following the notation of Theorem \ref{t5}, we denote
the resulting $48 \times 84$ matrix by $A_2$ (as in (\ref{a1a2})).

The rows of the $64 \times 84$ incidence matrix of any affine resolvable 
2-$(64,16,5)$ design $\cal{E}$ with a good block $B$ such that the residual design
${\cal{E}}_B$ coincides with ${\cal{D}}_B$=$(AG_{2}(3,4))_B$, are codewords of weight 21 in a binary 
linear code spanned by the rows of $A_2$ and one additional 
row $y=(y_1,\cdots,y_{84})$, where $y_{84}=1$ and the remaining 20 nonzero positions
of $y$ are labeled by the blocks of 5 parallel classes from the resolution $R$ of ${\cal{D}}''$.
In other words, $y$ is a row of the $16 \times 84$ matrix $A_1$ (cf. (\ref{a1a2})), labeled by
a point of $B$. 
Without loss of generality, we can fix one of the five parallel 
classes associated with the support of $y$. A computer check shows that among the
\[ { 19 \choose 4 }=3876 \]
choices for the remaining 4 parallel classes associated with the support of $y$,
only 16 lead to a code of length 84 and dimension 16 that contains sufficiently many
codewords of weight 21 to form the incidence matrix of a 2-$(64,16,5)$ design,
and each of these 16 codes does contain the incidence matrix of 
an affine resolvable 2-$(64,16,5)$  design. Further comparison shows that the set
 of 16 designs obtained
from the 16 codes contain two isomorphism classes of designs:  four designs
are isomorphic to $\cal{D}$=$AG_{2}(3,4)$ and 12 designs are isomorphic
to an affine resolvable  2-$(64,16,5)$ design ${\cal{E}}_1$ having 
full automorphism
group of order 92,160. The design ${\cal{E}}_1$ is isomorphic
to the affine resolvable design with these parameters 
found by Harada, Lam, and the author as the design
supported by minimum weight codewords of symmetric net No. 20
in \cite{HLT}.

The automorphism group of  ${\cal{E}}_1$  partitions its blocks into
three orbits, of length 1, 3, and 80 respectively.
The blocks from the orbit of length 80 are not good (in the sense of Definition
\ref{good}),  while the blocks from the orbits of length  one and three are good.

The residual design of  ${\cal{E}}_1$ with respect to the fixed block (orbit of length 1),
is isomorphic to a residual design of $AG_{2}(3,4)$, hence it admits two nonisomorphic
 linear embeddings: 
one in  $AG_{2}(3,4)$, and another in  ${\cal{E}}_1$.
 
A residual design of  ${\cal{E}}_1$ with respect to a good block $B'$ from the orbit
of length three is not isomorphic to  a residual design of $AG_{2}(3,4)$.
Its subdesign ${{\cal{E}}_1}''$ consisting of all blocks of size 12 
also has 32 resolutions, which split into three orbits of lengths 2, 10, and 20.
The weight distribution of the binary linear $[80,15]$ codes spanned by 
the incidence matrix of
${{\cal{E}}_1}''$ is given in Table \ref{tab2}.

\begin{tab}
\label{tab2}
\end{tab}
\begin{tabular}{|r|r|r|r|r|r|r|r|r|r|r|r|r|r|r|r|}
\hline
$i$ & 20 & 30 & 32 &  36 & 38 & 40 &  44 & 46 & 48 & 52 & 64 \\
\hline
$A_i$ & 48 &  1024 & 610 &  6240 & 10240 & 2880 & 5760 & 5120 & 600 & 240 & 5\\
\hline
\end{tabular}

\vspace{5mm}

We note that although the weight distribution in Table \ref{tab2} is different from that in Table \ref{tab1},
the number of codewords of weight 32 is again 610. Among these 610 codewords, there are 130 codewords
whose support is labeled by the blocks of a union of parallel classes of a resolution  form the orbit of
length 2. The number of codewords of weight 32 whose support is labeled by a union of parallel classes 
of a resolution from the orbit of length 10 or 20 is smaller than
\[ {16 \choose 2} = 120. \]
Following the procedure that we used for finding linear embeddings
of the residual design of $AG_{2}(3,4)$, an examination of the possible
choices for an additional row $y$ of weight 21 corresponding to a resolution of length 2,
establishes that in addition to ${\cal{E}}_1$, a residual design of  ${\cal{E}}_1$ with respect 
to a good block $B'$ from the orbit of length three can be embedded in an affine resolvable 2-$(64,16,5)$
design  ${\cal{E}}_2$ with full automorphism group of order 368,640, hence this design
is not isomorphic to $AG_{2}(3,4)$ or  ${\cal{E}}_1$. 

The design ${\cal{E}}_2$
is isomorphic to the design with the same parameters arising from net No. 36 in \cite{HLT},
as well as to the design arising from a special spread of lines in $PG(5,2)$
found by Mavron, McDonough and the author \cite{MMT}.

The automorphism group of  ${\cal{E}}_2$  partitions its blocks into
two orbits, of length 80 and 4, respectively. The blocks from the long orbit are not good,
while the blocks from the orbit of length 4 are good.
A residual design of ${\cal{E}}_2$ with respect to a good block
is isomorphic to a residual design of ${\cal{E}}_1$ (with respect to a block
from the orbit of length three), hence it has two nonisomorphic linear embeddings, in
${\cal{E}}_1$ and ${\cal{E}}_2$.

The next theorem  summarizes these results.
\begin{Thm}
\label{ag34}
(i) A residual design of $AG_{2}(3,4)$ admits exactly two
nonisomorphic linear embeddings over $GF(2)$: one in $AG_{2}(3,4)$, and 
a second one in an affine resolvable 2-$(64,16,5)$ design ${\cal{E}}_1$
of 2-rank 16 and having full automorphism group of order 92,160, 
 which is isomorphic to the affine resolvable design with 
these parameters arising from net No. 20 in \cite{HLT}. \\
(ii) The design  ${\cal{E}}_1$ has two types of good blocks.
 A residual design of  ${\cal{E}}_1$ with respect to a good block of the first type
 is linearly embeddable over $GF(2)$ in either   ${\cal{E}}_1$ or in $AG_{2}(3,4)$.
A residual design of  ${\cal{E}}_1$ with respect to a good block of the second type
is linearly embeddable over $GF(2)$ in 
either   ${\cal{E}}_1$, or in  an affine resolvable 2-$(64,16,5)$ design ${\cal{E}}_2$
of 2-rank 16 and having
full automorphism group of order  $368,640$, which is isomorphic to
the design with the same parameters arising from net No. 36 in \cite{HLT},
as well as to the design arising from a special spread in  $PG(5,2)$
 \cite{MMT}.\\
(iii)  A residual design of  ${\cal{E}}_2$ with respect to a good block
 is linearly embeddable 
over $GF(2)$ in either  ${\cal{E}}_2$ or ${\cal{E}}_1$.
\end{Thm} 

\begin{Rem}
{\rm
Hamada's conjecture \cite{Ha}, \cite{Ha81} states that the $p$-rank
of a design $\cal{D}$ having the same parameters as $PG_{d}(n,p^t)$ or
$AG_{d}(n,p^t)$, is greater than or equal to the $p$-rank of
 $PG_{d}(n,p^t)$ or $AG_{d}(n,p^t)$ respectively, with equality if and
only if $\cal{D}$ is isomorphic to  $PG_{d}(n,p^t)$ or
$AG_{d}(n,p^t)$. The $p$-ranks of  $PG_{d}(n,p^t)$ and
$AG_{d}(n,p^t)$, where $p$  a prime, $t\ge 1$,  and $n > d \ge 1$,
were computed by Hamada \cite{Ha}. 
Two affine resolvable 2-$(64,16,5)$ designs,
having the same parameters and the same 2-rank 16 as $AG_2(3,4)$, 
isomorphic to
 ${\cal{E}}_1$ and  ${\cal{E}}_2$ respectively,
were found originally by Harada, Lam and the author \cite{HLT}
as designs supported by  minimum weight codewords
in binary linear codes of length 64 and dimension 16, spanned by the
$64  \times 64$ incidence matrices of resolvable 1-$(64,16,16)$ designs 
whose dual designs are also resolvable, (or symmetric $(4,4)$-nets, 
in the terminology of \cite{HLT} and \cite{BJL}). These two designs are the 
only known
counter-examples to the "only-if" part of Hamada's conjecture
over a field of non-prime order \cite{CJT},  \cite{JT09}, \cite{T86}.
}
\end{Rem}

In view of Theorem \ref{ag34}, it will be interesting to know if
a residual design of  $AG_{n-1}(n,4)$ 
admits more than one nonisomorphic linear embeddings over $GF(2)$
in affine resolvable designs having the parameters of $AG_{n-1}(n,4)$,
for any value of $n$ greater than 3.

If $n=4$, 
 the substructure  ${\cal{D}}''$ of a 
 a residual design ${\cal{D}}_B$
with respect to a block $B$ of $AG_{3}(4,4)$
 consisting of all blocks of size 48,
 has  168 parallel classes and 2,097,152 resolutions.

By Hamada's formula \cite{Ha}, the 2-rank of  $AG_{3}(4,4)$ is 25.
The binary linear code of length 336 and dimension 24 spanned by the rows
of the incidence matrix $M''$ of ${\cal{D}}''$  
contains 10,290 codewords of weight 128, of which 2,226 codewords 
have supports being a union of 32 parallel classes of the resolution of
 ${\cal{D}}''$ defined by the blocks of the derived design
${\cal{D}}^B$. Since
\[ 2,226 > { 64 \choose 2} = 2016, \]
a residual design of  $AG_{3}(4,4)$ satisfies 
the condition of Theorem \ref{t5} with strict inequality,
and may have more than one linear embedding
over $GF(2)$.
However, finding all such embeddings by using the procedure applied 
to a residual design of $AG_{2}(3,4)$ seems to be  computationally infeasible.
One possible way to reduce the computations and make the problem problem manageable 
is by restricting the search to linear  embeddings which are invariant under 
a sufficiently large subgroup of the automorphism group of a residual design
of $AG_{3}(4,4)$.

\section{Affine resolvable designs as residual designs}

A 2-$(v,k,\lambda)$ design with replication number $r$
is called {\it quasi-residual} if it has the  parameters of
a residual design of a symmetric 2-$(v',k',\lambda')$ design with
\[ v' =v+r, \ k' = r, \ \lambda' = \lambda, \]
or equivalently, if
\[ r=k+\lambda. \]
By this definition, every affine resolvable
2-$(q^2\mu, q\mu, \frac{q\mu -1}{q-1})$ design is quasi-residual, 
and the parameters of a corresponding symmetric 2-design are
\begin{equation}
\label{ss}
 v' =\frac{q^3\mu-1}{q-1}, \ k'=\frac{q^2\mu-1}{q-1}, \ \lambda'=
\frac{q\mu-1}{q-1}.
\end{equation}
An affine  resolvable
2-$(q^2\mu, qm, \frac{q\mu -1}{q-1})$ design $\cal{D}$
 is the residual design 
of a symmetric 2-design ${\cal{D}}_1$ with parameters (\ref{ss}) 
if and only if 
there exists a symmetric 2-$(\frac{q^2\mu-1}{q-1},\frac{q\mu-1}{q-1},
\frac{\mu-1}{q-1})$ design ${\cal{D}}_0$ (S. S. Shrikhande \cite{SSS}, 
\cite[Corollary 5.4.9, page 178]{IS}).
A symmetric design with parameters (\ref{ss}) having $\cal{D}$ 
as a residual design is obtained by adding a block $B$ consisting of
$(q^2\mu-1)/(q-1)$ new points, being the points of
${\cal{D}}_0$,  choosing a bijection $\phi$ between the 
parallel classes of  $\cal{D}$ and the blocks
of ${\cal{D}}_0$,
 and  extending the $q$ blocks from any parallel class $P$
with the points of the block $\phi(P)$ of  ${\cal{D}}_0$.
By this construction, the derived design of  ${\cal{D}}_1$
with respect to the block $B$ is a
 2-$(\frac{q^2\mu-1}{q-1},\frac{q\mu-1}{q-1}, q\frac{\mu-1}{q-1})$ design, 
being a union of $q$ identical copies of
${\cal{D}}_0$. In the terminology of \cite{IS}, ${\cal{D}}_0$ is a
{\it normal} subdesign of ${\cal{D}}_1$. 

\begin{Def}
\label{def3}
{\rm
We call a block $B$ of a symmetric 2-design  ${\cal{D}}_1$
with parameters (\ref{ss})  {\it normal} (cf. \cite[XII.5]{BJL})
if the derived design of  ${\cal{D}}_1$ with respect to $B$ is
a union  of $q$ identical copies of a symmetric 
2-$(\frac{q^2\mu-1}{q-1},\frac{q\mu-1}{q-1}, q\frac{\mu-1}{q-1})$ design.
}
\end{Def}
 
An example of a symmetric design with all blocks being normal
is $PG_{n-1}(n,q)$, $n\ge 3$. In this case, any derived design 
is a union of $q$ identical copies of a design isomorphic to 
$PG_{n-2}(n-1,q)$, and any residual design is isomorphic to $AG_{n-1}(n,q)$.

\begin{Lem}
\label{lem1}
Let $q=p^t$, where $p$ is a prime, and let
${\cal D}$ be an affine resolvable 2-$(q^n, q^{n-1}, (q^{n-1}-1)/(q-1))$
design with $n\ge 2$ and a point by block incidence matrix $A$.
The design ${\cal D}$ is linearly embeddable over $GF(p)$
as a residual design in a symmetric
2-$((q^{n+1}-1)/(q-1), (q^n -1)/(q-1), (q^{n-1}/(q-1))$ design
 ${\cal{D}}_1$ if and only if the rows of a point by block incidence matrix
$A_1$ of ${\cal{D}}_1$  are codewords in the linear code over $GF(p)$ being 
the row space of the matrix (\ref{mat}).
\begin{equation}
\label{mat}
 \left(\begin{array}{cccc} & & & 0\\
&&& \cdot \\
& A & & \cdot \\
&&& \cdot \\
1 & \ldots & 1 & 1
\end{array}
\right). \end{equation}
\end{Lem}

{\bf Proof}. Since every column of $A_1$ contains
\[   \frac{q^n -1}{q-1}=q^{n-1}+\ldots + q +1 \equiv 1 \ {\pmod p}, \] 
the rows space of $A_1$ contains
the all-one vector $(1,1,\ldots 1)$, which is also the last row of (\ref{mat}).
$\Box$

We now apply Lemma \ref{lem1} to the affine resolvable 2-$(64,16,5)$ designs
 ${\cal{E}}_1$ and  ${\cal{E}}_2$ of 2-rank 16 discussed in Section \ref{sec5}.
If $A$ is the incidence matrix of   ${\cal{E}}_1$ or  ${\cal{E}}_2$,
then the row space of (\ref{mat}) over $GF(2)$ is a binary linear code of length 85 
and dimension 17 which contains exactly 69 codewords (thus, less that 85)
of weight 21. This implies the following.

\begin{Thm}
\label{6-3}
The affine resolvable  2-$(64,16,5)$ designs 
${\cal{E}}_1$ and  ${\cal{E}}_2$ do not admit any linear embedding
over $GF(2)$ in a symmetric 2-$(85,21,5)$ design.
Consequently, any symmetric 2-$(85,21,5)$ design having ${\cal{E}}_1$
or ${\cal{E}}_2$ as a residual design, must have 2-rank greater that 17,
which is the 2-rank of $PG_2(3,4)$.
\end{Thm}

\begin{Thm}
\label{t}
Let $q=p^t$, where $p$ is a prime. The classical affine resolvable design
$AG_{n-1}(n,2)$, $n\ge 2$, admits a unique linear embedding over $GF(p)$ in
$PG_{n-1}(n,q)$.
\end{Thm}

{\bf Proof}. The statement follows from Lemma \ref{lem1}, Theorem \ref{t4},
and the fact that the $p$-ary code spanned by an incidence matrix of
$PG_{n-1}(n,p^t)$ is of minimum weight $(q^n -1)/(q-1)$, and
by the restricted Johnson bound \cite[2.3.1]{HP}, \cite[2.4.2]{T88} , this code
cannot contain more than $(q^{n+1} -1)/(q-1)$ $(0,1)$-codewords of minimum 
weight. $\Box$

The following statement is an analogue of Theorem \ref{t5}
and gives a necessary condition for linear embeddability of affine
resolvable designs.

\begin{Thm}
\label{t-af}
Let $\cal{D}$  be an affine resolvable 2-$(q^n, q^{n-1}, (q^{n-1}-1)/(q-1))$ design,  $n\ge 2$,
where $q=p^t \ge 4$, and $p$ is a prime.
If $\cal{D}$ is linearly embeddable
over $GF(p)$ as a residual design in a symmetric 2-$((q^{n+1}-1)/(q-1), (q^{n}-1)/(q-1), (q^{n-1}-1)/(q-1))$
design with respect to a normal block, then the linear code over $GF(p)$ spanned by the rows
of the point-by-block incidence matrix of $\cal{D}$ contains at least
\[ (p-1){ \frac{q^n -1}{q-1} \choose 2} \]
codewords of weight $2q^{q-1}$ whose supports are unions of parallel classes of $\cal{D}$.
\end{Thm}

The proof is similar to that of Theorem \ref{t5}, so we omit it.\\

A quick computer check shows that the binary code of length 84 and dimension 16
 spanned by the incidence matrix of $AG_{2}(3,4)$
contains exactly
\[ 210 = { 21 \choose 2} \]
codewords of weight 32 whose supports are unions of parallel classes,
while the binary codes spanned by the incidence matrices of the affine resolvable 2-$(64,16,5)$ designs
${\cal{E}}_1$ and ${\cal{E}}_2$ contain only 130 such codewords. This provides another proof
of Theorem \ref{6-3}.

\section{ Acknowledgments} 

The author thanks Dieter Jungnickel for reading a preliminary version 
of this paper and making several useful remarks.
This research was supported by NSA Grant H98230-16-1-0011.

\end{document}